\documentclass{article}
\usepackage[final]{neurips_2025}

\PassOptionsToPackage{numbers, compress}{natbib}

\usepackage[utf8]{inputenc} % allow utf-8 input
\usepackage[T1]{fontenc}    % use 8-bit T1 fonts
\usepackage{hyperref}       % hyperlinks
\usepackage{url}            % simple URL typesetting
\usepackage{booktabs}       % professional-quality tables
\usepackage{amsfonts}       % blackboard math symbols
\usepackage{nicefrac}       % compact symbols for 1/2, etc.
\usepackage{microtype}      % microtypography
\usepackage{xcolor}         % colors
\usepackage{amsmath}
\usepackage[inline]{enumitem}
\usepackage{changepage}
\usepackage{todonotes}
\usepackage{tikz}
\usepackage{multirow}
\usetikzlibrary{arrows.meta, positioning, fit, shapes.geometric, calc}
\usepackage{graphicx}
\usepackage{wrapfig}
\usepackage{subcaption}
\usepackage[font=footnotesize]{caption}

\title{Learning to Optimize at Scale: A Benders Decomposition-TransfORmers Framework for Stochastic Combinatorial Optimization}

\author{%
  Seung Jin Choi \\
  Department of Industrial and Systems Engineering \\
  Virginia Tech \\
  Blacksburg, VA 24061\\
  \texttt{seungj9@vt.edu} \\
  \And
  Kimiya Jozani \\
  Department of Industrial and Systems Engineering \\
  Virginia Tech \\
  Blacksburg, VA 24061\\
  \texttt{kimiya@vt.edu} \\
  \AND
  Josh Cooper \\
  Department of Industrial and Systems Engineering \\
  Virginia Tech \\
  Blacksburg, VA 24061\\
  \texttt{joshfcooper@vt.edu} \\
    \AND
  Esra Buyuktahtakin Toy \\
  Department of Industrial and Systems Engineering \\
  Virginia Tech \\
  Blacksburg, VA 24061\\
  \texttt{esratoy@vt.edu} \\
}

\begin{document}
% \workshoptitle{MLxOR: Mathematical Foundations and Operational Integration of Machine Learning for Uncertainty-Aware Decision-Making}

\maketitle
\renewcommand\thefootnote{}\footnotetext{ \hspace{-2em} Presented at the NeurIPS 2025 Workshop on Machine Learning and Operations Research (MLxOR). Non-archival workshop paper.}

\addtocounter{footnote}{-1}

\begin{abstract}
We propose a learning-augmented Benders decomposition framework to solve large-scale two-stage stochastic mixed-integer programs. We focus on the two-stage stochastic capacitated lot-sizing problem (TSSCLSP) under demand uncertainty. Our method accelerates the convergence of the decomposition by using a pre-trained TransfORmer model to rapidly generate high-quality approximate solutions for the scenario subproblems. This hybrid strategy uses the TransfORmer's predictions to generate strong optimality and feasibility cuts, effectively guiding the Benders master problem. Our framework includes a novel expandable generation mechanism, allowing a model trained on a fixed horizon to solve instances of arbitrary length. For the test set considered, our method solves instances up to T=270, a scale previously intractable for this approach, while maintaining zero infeasibility in the generated subproblem solutions. This demonstrates the potential of TransfORmers as powerful surrogate solvers embedded within classical decomposition algorithms. \looseness-1
\end{abstract}

\section{Introduction}
The stochastic capacitated lot sizing (SCLSP) is a classic NP-hard problem that focuses on optimizing production planning in amidst uncertain demand and capacity limits. Traditional approaches, such as stochastic dynamic programming and decomposition techniques, face difficulties when applied at larger scales \citep{Birge2011,Shapiro2009,Sox1997,Askin1981}. Classical techniques, including Benders' decomposition \citep{Benders1962, VanSlyke1969}, Lagrangian relaxation \citep{sen2005algorithms, wolsey2020integer}, hedging algorithms \citep{Rockafellar1991, Watson2011, Haugen2001}, Dantzig-Wolfe decomposition \citep{Dantzig1960}, sample average approximation \citep{Shapiro2009, Kleywegt2002}, and scenario reduction \citep{Heitsch2003} have been employed to tackle the complexity of SCLSPs, but still struggle in complex problems as the problem scales.

Recently, the integration of deep learning (DL) models, such as TransfORmers, in optimization problems has gained momentum to approximate feasible regions or optimal policies in combinatorial settings \citep{Bengio2021, Khalil2017, bushaj2024k, bushaj2023simulation, cooper2024toward}. DL has shown significant potential in large-scale and real-time planning scenarios, such as SCLSP \citep{VanHezewijk2023, Temizoz2025, Felizardo2024, yilmaz2023learning,Gong2024, Karcher2025}. Transformers \citep{vaswani2017attention} have recently emerged as a powerful tool for sequential decision-making, showing success in large-scale settings \citep{yilmaz2023learning, lim2021time, cooper2024toward, yilmaz2024expandable, yilmaz2024deep, yilmaz2024non}. Their self-attention mechanism enables efficient parallel processing and modeling of long-range dependencies, making them ideal for high-dimensional, dynamic problems, such as the SCLSP. 

%%% MODIFIED SECTION: Aligned with Benders and expandable generation %%%
To address the scalability challenges of large-scale stochastic optimization and building on the emerging paradigm of ML-OR~\citep{cooper2024toward, yilmaz2024expandable, yilmaz2024non}, we propose a novel \textbf{TransfORmer-accelerated Benders decomposition} framework to solve the two-stage stochastic capacitated lot-sizing problem (TSSCLSP). Our method addresses the primary bottleneck of Benders decomposition---the iterative solving of numerous scenario subproblems---by replacing the exact subproblem solver with a TransfORmer-based surrogate. This model rapidly predicts near-optimal subproblem solutions, which are then used to generate strong optimality and feasibility cuts to guide the master problem. Critically, we introduce an \textbf{expandable, sliding-window generation technique} that enables a Transformer trained on a fixed-length horizon (e.g., T=90) to produce complete solutions for significantly longer, unseen horizons (e.g., T=270). This framework advances the state-of-the-art by coupling exact mathematical decomposition with scalable, learning-based inference, bridging the gap between optimization and machine learning. Although this paper focuses on the TSSCLSP, the approach is general and can be adapted to a wide range of dynamic two-stage stochastic optimization problems. \looseness-1

\section{Mathematical Formulation}\label{sec:mathmod}
The mathematical model of TSSCLSP is presented in~\eqref{eq1}--\eqref{eq5}~\citep{buyuktahtakin2023scenario}.
Let $\mathcal{T} = \{1,\ldots,T\}$ denote the set of time periods. Let $\omega \in \Omega$ index the demand scenarios. The deterministic parameters are: $f_t$ (fixed production cost), $g_t$ (unit production cost), and $\kappa_t$ (production capacity) in period $t$. The stochastic parameters are: $d_t^\omega$ (demand in scenario $\omega$), $h_t$ (inventory holding cost) and $p^\omega$ (probability of scenario $\omega$). The decision variables are: $Y_t \in \{0,1\}$, a first-stage binary variable that indicates whether production is initiated; $X_t^\omega \geq 0$ and $S_t^\omega \geq 0$ are the second-stage continuous variables for quantity produced and inventory in period $t$ under scenario $\omega$. \looseness-1
 
\begin{minipage}{0.56\textwidth}
{\footnotesize\begin{flalign}
  \min \quad& \sum_{t\in\mathcal{T}} f_t Y_t+ \sum_{\omega\in\Omega}p^\omega\left(\sum_{t\in\mathcal{T}} g_t X_t^\omega + h_t S_t^\omega \right) & \label{eq1} 
  \end{flalign}
  \vspace{-1em}
  \begin{flalign}
  \text{s.t.} \quad & S_{t-1}^\omega + X_t^\omega - d_t^\omega = S_t^\omega & \forall \omega \in \Omega , \forall t \in \mathcal{T} \label{eq2} \\
  & X_t^\omega \leq \kappa_t Y_t & \forall t\in \mathcal{T}, \forall \omega \in \Omega \label{eq3} \\
  & X_t^\omega, S_t^\omega \geq 0 & \forall \omega \in \Omega , \forall t \in \mathcal{T} \label{eq4} \\
  & Y_t \in \{0, 1\} & \forall t \in \mathcal{T}. \label{eq5}
\end{flalign}}
\end{minipage}
\hfill
\begin{minipage}{0.4\textwidth}
The objective function \eqref{eq1} minimizes the total first-stage fixed setup costs and the expected second-stage production and inventory holding costs. Constraint \eqref{eq2} enforces inventory balance. Constraint \eqref{eq3} links production to the setup decision, subject to capacity $\kappa_t$. Constraints \eqref{eq4} and \eqref{eq5} define the variable domains.
\end{minipage}

%%% MODIFIED SECTION: Replaced dynamic decomposition with Benders framework %%%
\section{Methodology}\label{sec:meth}
The extensive form in \eqref{eq1}--\eqref{eq5} is $NP$-hard and intractable for large-scale instances. We therefore employ a Benders decomposition, which we accelerate leveraging a TransfORmer model that solves scenario sub-problems.

\subsection{Benders Decomposition Framework}
The problem is decomposed into a Master Problem (MP) that determines the first-stage setup decisions ($Y_t$) and a set of independent, scenario-specific Subproblems (SP) that evaluate the consequences of those decisions.
The MP minimizes setup costs plus an approximation of the expected recourse cost, represented by variables $\eta_\omega$.
\begin{align}
    \textbf{(MP)} \quad \min \quad & \sum_{t\in\mathcal{T}} f_t Y_t + \sum_{\omega\in\Omega} p^\omega \eta_\omega \\
    \text{s.t.} \quad & Y_t \in \{0,1\} \quad \forall t \in \mathcal{T} \\
   % & \eta_\omega \ge [\text{Optimality and Feasibility Cuts}] \quad \forall \omega \in \Omega
       & \eta_\omega \ge \theta_\omega^k + (\pi_\omega^k)^\top (Y - \hat{Y}^k), 
      \quad \forall \omega \in \Omega, \; \forall k \in \mathcal{K}_\omega. 
      \label{eq:MP_cuts}
\end{align}
Constraint \eqref{eq:MP_cuts} represents the family of Benders optimality and feasibility cuts that iteratively approximate the expected recourse function, where $(\pi_\omega^k, \theta_\omega^k)$ are dual solutions obtained from the $k$th subproblem of scenario $\omega$. For a fixed master solution $\hat{Y}$, each SP reduces to a Linear Program (LP); its dual is solved to generate these cuts, which are then added to the MP to iteratively refine the solution.

\subsection{Learning-Augmented Subproblem Solving with TransfORmers}
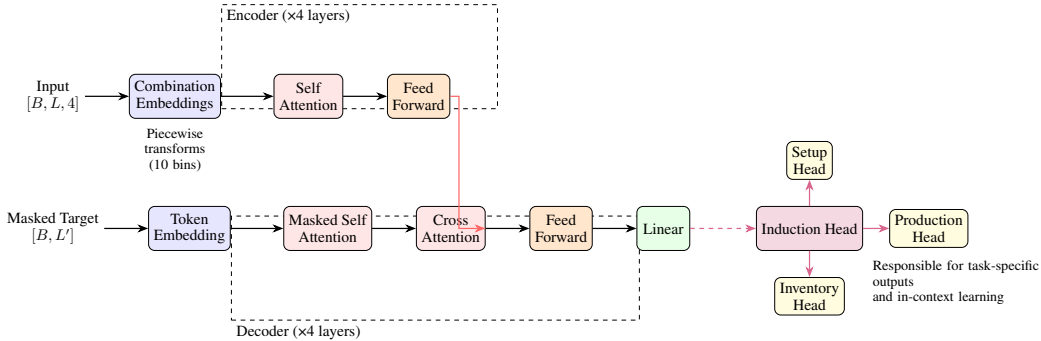
\begin{figure}[ht]
  \centering
  \resizebox{\linewidth}{!}{%
    \begin{tikzpicture}[%
        node distance = 0.4cm and 1.0cm,
        box/.style   = {draw, minimum width = 1.2cm, minimum height = 1cm,
                        align = center, rounded corners, font=\normalsize},
        embedding/.style = {box, fill = blue!10},
        attention/.style = {box, fill = red!10},
        ff/.style        = {box, fill = orange!20},
        linear/.style    = {box, fill = green!10},
        induction/.style = {box, fill = purple!15},
        task/.style      = {box, fill = yellow!15, minimum width = 1.0cm, minimum height =.8cm},
        arrow/.style     = {-{Stealth[length=2.5mm]}, thick},
        data/.style      = {font = \normalsize, align = center}
    ]
% Title
% Input
\node[data] at (-15,0) (input) {Input\\$[B, L, 4]$};
% Embedding
\node[embedding, right=of input] (feature_embed) {Combination\\Embeddings};
% Outer box for Encoder Block
\node[draw, dashed, fit={($(feature_embed.east)+(0.3,0)$) ($(feature_embed.east)+(6,1.8)$)}, inner sep=8pt, label={[anchor=north west, font=\normalsize]north west:Encoder (×4 layers)}] (enc_block) {};
% Encoder components
\node[attention, right=1.2cm of feature_embed] (enc_attn) {Self\\Attention};
\node[ff, right=of enc_attn] (enc_ff) {Feed\\Forward};
% Target Input
\node[data] at (-15,-3) (tgt_in) {Masked Target\\$[B, L']$};
\node[embedding, right=of tgt_in] (tgt_embed) {Token\\Embedding};
% Outer box for Decoder Block
\node[draw, dashed, fit={($(tgt_embed.east)+(0.3,0)$) ($(tgt_embed.east)+(9,-1.8)$)}, inner sep=8pt, label={[anchor=north west, font=\normalsize]south west:Decoder (×4 layers)}] (dec_block) {};
% Decoder components
\node[attention, right=1.2cm of tgt_embed] (dec_self_attn) {Masked Self\\Attention};
\node[attention, right=of dec_self_attn] (dec_cross_attn) {Cross\\Attention};
\node[ff, right=of dec_cross_attn] (dec_ff) {Feed\\Forward};
% Output
\node[linear, right=of dec_ff] (output) {Linear};
% \node[data, right=of output] (final) {Output\\$[B, L', V]$}; % Removed to declutter

% NEW ELEMENT: Induction Head
\node[induction, right=1.5cm of output] (induction_head) {Induction Head};

% NEW ELEMENT: Task-specific outputs
\node[task, above=0.6cm of induction_head] (setup_head) {Setup\\Head};
\node[task, right=0.6cm of induction_head] (production_head) {Production\\Head};
\node[task, below=0.6cm of induction_head] (inventory_head) {Inventory\\Head};

% Connections - Encoder path
\draw[arrow] (input) -- (feature_embed);
\draw[arrow] (feature_embed) -- (enc_attn);
\draw[arrow] (enc_attn) -- (enc_ff);
% Connections - Decoder path
\draw[arrow] (tgt_in) -- (tgt_embed);
\draw[arrow] (tgt_embed) -- (dec_self_attn);
\draw[arrow] (dec_self_attn) -- (dec_cross_attn);
\draw[arrow] (dec_cross_attn) -- (dec_ff);
\draw[arrow] (dec_ff) -- (output);
% Cross-attention connection
\draw[arrow, red!60] (enc_ff) -- ($(enc_ff.east)+(0.2,0)$) |- (dec_cross_attn);
% NEW CONNECTIONS: Induction head
\draw[arrow, dashed, purple!60] (output) -- (induction_head);
\draw[arrow, purple!60] (induction_head) -- (setup_head);
\draw[arrow, purple!60] (induction_head) -- (production_head);
\draw[arrow, purple!60] (induction_head) -- (inventory_head);
\node[below=0.1cm of feature_embed, align=center, text width=1.8cm, font=\small] 
    {Piecewise\\transforms\\(10 bins)};
\node[below right=0.1cm and 0.1cm of induction_head, align=left, text width=4cm, font=\small] 
    {Responsible for task-specific outputs\\and in-context learning};
\end{tikzpicture}
  }
   \caption{TransfORmer architecture used to solve subproblems in the decomposition framework.}
     \label{fig:transformer_arch}
\end{figure}

The bottleneck in Benders decomposition is the repeated solving of SPs. We accelerate this by using a TransfORmer model as a surrogate learning-to-optimize solver. The process is illustrated in Figure \ref{fig:transformer_arch}. Within each Benders worker:
\begin{enumerate*}[label=(\roman*)]
    \item A partial solution from the MP is encoded into a masked input tensor.
    \item The TransfORmer performs a single forward pass to complete the solution, predicting the binary $Y_t$ variables for the full horizon of a specific scenario.
    \item With all $Y_t$ variables fixed, the SP reduces to a fast LP, which is solved to obtain exact costs and duals.
    \item These duals are used to construct valid Benders cuts, which are sent back to the master.
\end{enumerate*}

\subsection{Expandable Generation for Long-Horizon Problems}

\begin{wraptable}{r}{0.48\textwidth}
\vspace{-1.4em}
\centering
\caption{\footnotesize TransfORmer hyperparameter configuration.}
\label{tab:hyperparams}
{\footnotesize
\begin{tabular}{lc}
\toprule
\textbf{Parameter} & \textbf{Value} \\
\midrule
Parameter Count & 6.8M \\
Source Vocabulary Size & 16,000 \\
Target Vocabulary Size & 5 \\
Embedding Dimensions & 256 \\
Training GPU Hours & <1 \\
\bottomrule
\end{tabular}}
\vspace{-1.4em}
\end{wraptable}
A key contribution is an expandable generation mechanism that allows a TransfORmer trained on a fixed horizon (e.g., $T=90$) to solve instances of arbitrary length (e.g., $T=270$). This is achieved via a sliding window method:
\begin{enumerate*}[label=(\roman*)]
    \item \textit{Initial Generation}: The model generates a solution for the first window (e.g., periods 0-89).
    \item \textit{Slide and Predict}: To generate the next chunk, the model's input uses the last portion of the previous output as historical context, concatenated with any known future master variables and `[MASK]` tokens.
    \item \textit{Iterative Completion}: This process repeats, sliding the window until a complete solution for the full length is constructed.
\end{enumerate*}

The model configuration is summarized in Table~\ref{tab:hyperparams}. To increase robustness, we apply 5-shot sampling per subproblem and select the best feasible prediction to construct cuts.

%%% MODIFIED SECTION: Replace your entire Section 4 with this one. %%%

%%% MODIFIED SECTION: Replace your entire Section 4 with this new version. %%%

\section{Results of Numerical Experiments}\label{sec:res}
We evaluated the proposed decomposition-TransfORmer framework (ML-Benders) on two-stage TSSCLSP instances, using the benchmark style of \citet{atamturk2004study} and  \citet{yilmaz2023learning}. The TransfORmer model was trained exclusively on deterministic capacitated lot-sizing problem instances with a fixed horizon of $T=90$. Our experiments are designed to test two key hypotheses: first, the viability and robustness of our expandable generation mechanism on long-horizon problems and second, the quality of the TransfORmer-based cuts at the model's native training horizon.

\subsection{Scalability Demonstration on Long-Horizon Instances (T=270)}
The primary test of our framework is its ability to scale to problem sizes far beyond its training data. We tested the model, trained only on $T=90$ instances, on a set of TSSCLSP instances with a horizon of $T=270$. This represents a challenging zero-shot generalization task.
The key result is the successful application of the expandable generation mechanism. The framework was able to consistently generate complete, feasible solutions for all scenarios in these long-horizon problems. As shown in Table~\ref{tab:long-horizon-results}, this was achieved with only a modest increase in the number of Benders iterations compared to the baseline task.

\begin{table}[htbp]
\centering
\caption{\centering Performance of ML-Benders on the expandable long-horizon task. \\The model was trained only on $T=90$ instances.}
\label{tab:long-horizon-results}
\footnotesize
\begin{tabular}{lcc}
\toprule
\textbf{Test Horizon} & \textbf{Avg. Iterations} & \textbf{Final Opt. Gap (\%)} \\
\midrule
$T=90$ (Baseline Training Horizon) & 51.6 & 3.55\% \\
$\mathbf{T=270}$ \textbf{(Expandable Test)} & \textbf{66.9} & \textbf{19.60\%} \\
\bottomrule
\end{tabular}
\end{table}

The final average optimality gap of 19.60\% indicates a trade-off in solution quality when generalizing to a horizon three times the training length. However, the ability to generate valid, structured solutions for a problem of this complexity and scale—without any direct training and without causing subproblem infeasibility—is a significant demonstration of the model's robust generalization. This validates our expandable framework as a viable method for tackling large-scale instances that would otherwise be intractable.

\subsection{Analysis of Cut Quality at the Training Horizon (T=90)}
To understand the foundational performance that enables this scalability, we performed a detailed analysis comparing ML-Benders with classical Benders in 50 instances at the native training horizon of $T=90$. These results highlight the superior quality of the cuts generated by solving sub-problems with TransfORmer.

\subsubsection{Performance under a 300-Second Time Limit}
In a time-constrained setting, the ML-Benders framework demonstrates clear superiority across all metrics, as summarized in Table~\ref{tab:results-summary-90}.

\begin{table}[htbp]
\centering
\caption{\centering Average performance comparison under a 300-second time limit \\($T=90$, 10 scenarios, 50 instances).}
\label{tab:results-summary-90}
\footnotesize
\begin{tabular}{lccc}
\toprule
\textbf{Method} & \textbf{Avg. Iterations} & \textbf{Avg. Time (s)} & \textbf{Final Opt. Gap (\%)} \\
\midrule
Classical Benders & 73.1 & 308.2 & 41.81\% \\
\textbf{ML-Benders} & \textbf{51.6} & \textbf{246.8} & \textbf{3.55\%} \\
\bottomrule
\end{tabular}
\end{table}

The key findings are:
\begin{enumerate*}[label=(\roman*)]
    \item \textit{Superior Solution Quality:} At its native horizon, ML-Benders achieves an average gap of only \textbf{3.55\%}, a \textbf{91.5\% reduction} compared to classical Benders (41.81\%).
    \item \textit{Enhanced Efficiency:} The method requires \textbf{29.4\% fewer iterations} and \textbf{19.9\% less time}, confirming that TransfORmer-generated cuts are highly effective in guiding the master problem.
\end{enumerate*}

\subsubsection{Cut-Budget Sensitivity (50-Cut Limit)}
To directly test the quality of the generated cuts, we performed a "stress test" by imposing a strict budget of only 50 cuts. The computational results are presented in Table~\ref{tab:cut-limit-results}.

\begin{table}[htbp]
\centering
\caption{Average optimality gap under a strict 50-cut limit.}
\label{tab:cut-limit-results}
\footnotesize
\begin{tabular}{lc}
\toprule
\textbf{Method} & \textbf{Final Opt. Gap (\%)} \\
\midrule
Classical Benders & 77.03\%* \\
\textbf{ML-Benders} & \textbf{2.80\%} \\
\bottomrule
\end{tabular}
\caption*{\scriptsize *Average calculated over 9 instances; one instance failed to find a feasible solution.}
\end{table}

With a severely restricted budget, classical Benders fails to converge, whereas ML-Benders remains robust, achieving a final gap of just \textbf{2.80\%}. This provides conclusive evidence that the TransfORmer-based surrogate generates dense, highly informative constraints. It is this foundational strength in generating high-quality cuts that allows the expandable mechanism to function effectively on much larger problems.

\section{Discussion and Future Directions}\label{Discussion and Future Directions}
This work integrates a TransfORmer-based predictor with a Benders-style decomposition, yielding a hybrid that couples the efficiency of deep learning with the rigor of mathematical optimization. The success of our framework, particularly its ability to solve $T=270$ instances without generating infeasible subproblems, is a significant step towards building trust in ML-augmented optimization solvers. The \textit{expandable generation mechanism} is a key contribution, proving that models trained on smaller, manageable problem sizes can generalize their learned policies to solve much larger, industrially-relevant instances. The framework is problem-agnostic and extensible to other two-stage settings. Future work will focus on completing the large-scale stochastic experiments, refining architectures, and exploring adaptive retraining to further improve cut quality and robustness.

{\small
\bibliographystyle{plainnat}
\bibliography{refs}
}

\end{document}